\documentclass[12pt,fleqn] {article}
\usepackage{amssymb,amsmath,amsthm}
\usepackage{graphicx}
\usepackage{a4wide}
\usepackage{url}
\usepackage[natbibapa]{apacite}
\usepackage{doi}

\usepackage{etoolbox}
\AtBeginEnvironment{APACrefURL}{\renewcommand{\url}[1]{}}
\renewcommand{\APACrefYearMonthDay}[3]{\BBOP{#1}\BBCP}
\usepackage{hyperref}
\newcommand{\R}{\mathbb{R}}

\newcommand{\bs}{\boldsymbol}

\newcommand{\adb}{\allowdisplaybreaks}

\newcommand{\inv}{\frac{1}}
\date{\today}
\author{ Thomas Fung$^{a,}$\footnote{Corresponding Author. Email address:  \href{mailto:thomas.fung@mq.edu.au}{thomas.fung@mq.edu.au} (Thomas Fung). Honorary Associate, University of Sydney.}\, and Eugene Seneta$^{b}$\\
{\small $^a$ Department of Mathematics and Statistics, Macquarie University, NSW 2109, Australia}\\
{\small $^b$ School of Mathematics and Statistics, University of Sydney, NSW 2006, Australia}\\
}
\title{Tail asymptotics for the bivariate equi-skew Variance-Gamma distribution}
\begin{document}

\renewcommand{\APACrefYearMonthDay}[3]{\APACrefYear{#1}}

\maketitle

\begin{abstract}
\noindent  We derive the asymptotic rate of decay to zero of the tail dependence of the bivariate skew Variance Gamma (VG) distribution under the equal-skewness condition, as an explicit  regularly varying function.   Our development is  in terms of a slightly more general bivariate  skew Generalized Hyperbolic (GH) distribution. Our initial reduction of the bivariate problem to a univariate one is motivated by our earlier study of tail dependence rate for the bivariate skew normal distribution.
\section{Introduction}
\end{abstract}
The coefficient of lower tail dependence of a random vector $\textbf{X} =
(X_1,X_2)^{\top}$ with marginal inverse distribution functions $F_1^{-1}$ and
$F_2^{-1}$ is defined as
\begin{equation}
\lambda_L = \lim_{u \rightarrow 0^+}\lambda_L(u), \quad \text{where} \quad \lambda_L(u) = P(X_1 
\leq F_1^{-1}(u) | X_2
\leq F_2^{-1}(u)).\label{defn:tail dependence}
\end{equation}

\noindent $\textbf{X}$ is said to have asymptotic lower  tail dependence if $\lambda_L$ exists and is positive. If $\lambda_L=0$, then $\textbf{X}$ is said to be
asymptotically independent in the lower tail. 

This quantity provides insight on
the tendency for the distribution to generate joint extreme event since it
measures the strength of dependence (or association) in the lower tail of a
bivariate distribution.  If the marginal distributions of these random variables are
continuous, then from (\ref{defn:tail dependence}), it follows that
$\lambda_L(u)$ can be expressed in terms of  the copula of
$\textbf{X}$, $C(u_1,u_2)$, as
\begin{equation*}
\lambda_L(u) = \frac{P(X_1\leq F_1^{-1}(u), X_2\leq F_2^{-1}(u))}{P(X_2 \leq F_2^{-1}(u))} =  \frac{C(u,u)}{u}.
\label{defn:tail dependence 3}
\end{equation*}

If $\lambda_L =0 $ in (\ref{defn:tail dependence}), that is, if  asymptotic lower tail independence obtains, the asymptotic rate of convergence to zero  as $ u \to 0^+$ of the copula $C(u,u)$  is tantamount to that of $\lambda_L(u)$ 
through the relation:
\begin{equation*}
C(u,u)= u\lambda_L(u). \label{copula} 
\end{equation*}

The central purpose of this paper is to provide  an analytic result on the asymptotic tail independence for the bivariate skew Variance-Gamma (VG) model. We will consider this problem in terms of the more general skew Generalized Hyperbolic (GH) distribution.


The (standardized) bivariate skew GH distribution, $GH(0,R,\bs{\theta}, p,a,b)$ is defined by its variance-mean mixing representation as 
\begin{equation}
\textbf{X} = \bs{\theta} W + \sqrt{W}\textbf{Z}   \quad \label{define: skew GH}
\end{equation}
where $\textbf{X} = (X_1, X_2)^{\top}, \bs{\theta}^{\top} = (\theta_1,\theta_2)$, and $ W \sim GIG(p,a,b)$ is independently distributed of $\textbf{Z} \sim N(0, R)$. Here $R = \left( \begin{smallmatrix} 1 & & \rho \\ \\ \rho & & 1\end{smallmatrix}\right)$, with $-1<\rho<1$. 

Recall that a random variable $W$ is said to have a (univariate) Generalised Inverse
Gaussian (GIG) distribution, denoted by $GIG(p,a,b)$, if it has density
\begin{eqnarray*}
f_{GIG}(w) &=& \frac{1}{2\overline{K}_{p}(a,b)}w^{p-1}\exp(-\frac{1}{2}(a^2w^{-1}+b^2w)), \quad  w>0; \label{density:GIG}
\\
\notag &=& 0, \quad \text{otherwise;} \end{eqnarray*}
where
{\adb
\begin{equation}
 \overline{K}_{p}(a,b) = \begin{cases}
   (\frac{a}{b})^{p}K_{p}(ab),
   &\text{ $p \in\R$, if  $a,b >0$;}\\
   b^{-2p}\Gamma(p)2^{p-1}, &\text{
   $p ,b>0$, if $a=0$;}\\
   a^{2p}\Gamma(-p)2^{-p-1}, &\text{
   $a>0$ and $p <0$, if  $b=0$,}
 \end{cases}\label{Kbar}
\end{equation}

Here $K_{p}(\omega), \ \omega >0, $ is the modified Bessel function of the second kind (\cite{BatemanManuscriptProject.1954}) with index $p \in \R$.

In  the VG special  case   $a=0$, $b = \sqrt{\frac{2}{\nu}}$, $p=\inv{2}$.  We proceed  
more generally by assuming $ b > 0 $ in this note in the $GIG(p,a,b)$ setting.

It was shown in \cite{Fung2011a} that when $\textbf{X}\sim GH(0,R,\bs{\theta}, p,a,b)$ with $b>0$, then $\textbf{X}$ is asymptotically independent in the lower tail; that is  $\lambda_L= 0.$  The proof in \cite{Hammerstein2016} for the VG can be adjusted to give this same conclusion.

Our specific focus in the sequel is to obtain a rate of convergence result of the form: 
\begin{equation} \label{convergence rate 1}
C(u,u) = u^{\tau}L(u) 
\end{equation} where $L(u)$ is a slowly varying function (SVF) as $ u \to 0^{+}$ and $\tau >1,$ when $\theta_1=\theta_2, = \theta$, say, so $X_i \sim \theta W + \sqrt{W}Z_i$, $i=1,2$, where $Z_i \sim N(0,1).$ That is, the distribution functions of the $X_i, i =1,2$ are the same: $F_1(u) = F_2(u), = F(u)$, say.  We call this assumption in a bivariate setting ``equi-skewness".

Our study thus parallels that of \cite{Fung2016}, who treat the bivariate skew normal distributed $\textbf{X}$, that is $\textbf{X}\sim SN_2({\bf \alpha},R)$ where in 
${\bf \alpha}=(\alpha_1, \alpha_2)^{\top}$, it is  assumed that  $\alpha_1 = \alpha_2, = \alpha$ say, 
so equi-skewness obtains.  

Both treatments depend  on the same initial device: that 
\begin{equation} \label{max} 
Z_{(2)} = \max (Z_1,Z_2)  \sim SN_1(\alpha), {\text{ where}}\, \,  \alpha = \sqrt{ \frac{1- \rho}{1+\rho}},
\end{equation} to reduce the bivariate problem to a univariate one. Our subsequent treatment is quite different, since the setting in
\cite{Fung2016} is just mean-mixing, and with a mixing distribution not encompassed by the GIG.

Clearly, since $\lambda_L(u)$ is a probability, the index $\tau$ in (\ref{convergence rate 1}) must satisfy $\tau \geq 1.$ We  note that \citet{Ledford1997}, \citet{Ramos2009},  \citet{Hashorva2010} and \citet{Hua2011} all classified the degree of tail-dependence to the value of $\tau$ in  (\ref{convergence rate 1}).  \citet{Hua2011} define $\tau$ in (\ref{convergence rate 1}) as the (lower) tail order of a copula. The tail order case $1 < \tau<2$ is considered as  intermediate tail dependence as it corresponds to  the copula having some level of positive dependence in the tail when $\lambda_L=0$. Thus when $\lambda_L(u)=C(u,u)/u=u^{\tau-1}L(u)$, $1< \tau <2$, there is some measure of positive association when $\lambda_L =0$, but the association is not as strong as when $\tau =1$, and $\lambda_L(u)= L(u) \to \lambda_L >0$, $u \to 0^+$, the case of asymptotic tail dependence.
 In our specific setting we shall find that $ 1 <\tau < \infty.$


For the case of $b=0$ in the $GIG(p, a, b)$ setting, $\textbf{X}$ can be asymptotically dependent in the lower tail. The limiting and rate of convergence results for this case were discussed in \cite{Fung2010} and \cite{Fung2014} respectively.

\section{The Reduction}
For our equi-skew setting of (\ref{define: skew GH}):
\begin{eqnarray*}
&& P(X_1\leq F_1^{-1}(u), \, X_2\leq F_2^{-1}(u)) \\
&=& P(X_1 \leq y, \, X_2 \leq y), \quad \text{where $y= F_1^{-1}(u) = F_2^{-1}(u)$};\\
&=& E_W(P(\theta W + \sqrt{W} Z_1\leq y, \, \theta W + \sqrt{W} Z_2\leq y))\\
&=& E_W\left[P\left(Z_1 \leq \frac{y-\theta W}{\sqrt{W}},\, Z_2 \leq \frac{y-\theta W}{\sqrt{W}}\right)\right]\\
&=& E_W\left[ P(Z_{(2)} \leq \frac{y-\theta W}{\sqrt{W}})\right],
\end{eqnarray*}
where $Z_{(2)} $ has a skew normal distribution with skew parameter $\alpha$ according to (\ref{max});{\adb
\begin{eqnarray*}
&=& \int^{\infty}_{0} P\left(Z_{(2)} \leq \frac{y-\theta w}{\sqrt{w}}\right)f_W(w)\,dw\\
&=& \int^{\infty}_0 P(\theta w+\sqrt{w}Z_{(2)} \leq y) f_W(w)\,dw\\
&=& P(X^* \leq y),
\end{eqnarray*}}
where $X^* = \theta W +\sqrt{W} Z_{(2)}$ and $X^*$ is defined by a variance-mean mixing of a skew normal representation. This type of distribution was considered in \cite{Arslan2014}, according to whose Proposition 1 the probability density of $X^*$ is given by:
\begin{equation} \label{fxstar}
f_{X^*}(x)= \frac{2e^{\theta x}}{\sqrt{2\pi}\, \overline{K}_{p}(a,b)} \overline{K}_{p-\inv{2}}\left((a^2+x^2)^{\inv{2}}, (\theta^2+b^2)^{\inv{2}}\right) P(Y \leq \alpha x) \end{equation}
where $Y\sim GH(0,1,\alpha\theta, p-\inv{2},(a^2+x^2)^{\inv{2}},(\theta^2+b^2)^{\inv{2}})$ i.e. $Y$ has a univariate 
GH distribution. 
Now, setting:  $\, \, \beta = (\alpha^2\theta^2+\theta^2+b^2)^{\inv{2}}, = (\theta^2(1+ \alpha^2)+b^2)^{\inv{2}}$:  
{\adb
\begin{eqnarray*}
\notag && P(Y\leq \alpha x) \\
\notag&=& \int^{\alpha x}_{-\infty} \frac{e^{\alpha \theta z}}{\sqrt{2\pi}\,\overline{K}_{p- \inv{2}}((a^2+x^2)^{\inv{2}}, (\theta^2+b^2)^{\inv{2}})} \overline{K}_{p-1}\left( (a^2+x^2+z^2)^{\inv{2}},\beta \right)\,dz,  \\
&=& \int^{\alpha x}_{-\infty}\frac{e^{\alpha \theta z}}{\sqrt{2\pi}\,\overline{K}_{p- \inv{2}}((a^2+x^2)^{\inv{2}}, (\theta^2+b^2)^{\inv{2}})}\left(\frac{(a^2+x^2+z^2)^{\inv{2}}}{\beta}\right)^{p-1} K_{p-1}\left( \beta (a^2+x^2+z^2)^{\inv{2}} \right)\, dz 
\end{eqnarray*}}
from (\ref{Kbar}). As $z\leq x\leq y$, so when $y\to -\infty$, $(a^2+x^2+z^2)^{\inv{2}} \to \infty$, we can use the asymptotic behaviour of the Bessel function (see \cite{Jorgensen1982}):
\begin{equation} \label{bessel}
K_{\nu}(y) =  \sqrt{\frac{\pi}{2 y}}e^{-y}\left(1 + O\left(\inv{y}\right)\right), \quad \text{as $y\to \infty$}
\end{equation}
and $P(Y\leq \alpha x)$ becomes  {\adb
\begin{align*}
= & \int^{\alpha x}_{-\infty}\frac{e^{\alpha \theta z}}{\sqrt{2\pi}\,\overline{K}_{p- \inv{2}}((a^2+x^2), (\theta^2+b^2)^{\inv{2}})}\left(\frac{(a^2+x^2+z^2)^{\inv{2}}}{\beta}\right)^{p-1}\\
& \quad \times  \sqrt{\frac{\pi}{2\beta(a^2+x^2+z^2)^{\inv{2}}}} e^{-\beta(a^2+ x^2+z^2)^{\inv{2}}}
\left(1+O\left(\inv{\sqrt{a^2+x^2+z^2}}\right)\right)\,dz\\
=& \int^{\infty}_{\alpha |x|} \frac{e^{- \alpha \theta z}}{\sqrt{2\pi}\,\overline{K}_{p-\inv{2}}((a^2+x^2)^{\inv{2}}, (\theta^2+b^2)^{\inv{2}})}\left(\frac{(a^2+x^2+z^2)^{\inv{2}}}{\beta}\right)^{p-1} \\
& \quad \times \sqrt{\frac{\pi}{2\beta(a^2+x^2+z^2)^{\inv{2}}}} e^{-\beta(a^2+x^2+z^2)^{\inv{2}}}\left(1+O\left(\inv{|x|}\right)\right)\,dz\\
=& \int^{\infty}_{\alpha } \frac{e^{- \alpha \theta |x|s}}{\sqrt{2\pi}\,\overline{K}_{p-\inv{2}}((a^2+x^2)^{\inv{2}}, (\theta^2+b^2)^{\inv{2}})}\left(|x| \frac{\left(1+\frac{a^2}{x^2}+s^2\right)^{\inv{2}}}{\beta}\right)^{p-1} \\
& \quad \times \sqrt{\frac{\pi}{2\beta|x|\left(1+\frac{a^2}{x^2}+s^2\right)^{\inv{2}}}} 
e^{-\beta |x| \left(1+\frac{a^2}{x^2}+s^2\right)^{\inv{2}}}|x|
\left(1+O\left(\inv{|x|}\right)\right)\,ds,
\end{align*}}
by letting $z = |x| s$;
\begin{align}
\notag =& \frac{|x|^{p-\inv{2}}}{2\,\overline{K}_{p- \inv{2}}((a^2+x^2)^{\inv{2}}, (\theta^2+b^2)^{\inv{2}}) \beta^{p-\inv{2}}} 
\left(1+O\left(\inv{|x|}\right)\right) \\
& \quad \times \int^{\infty}_{\alpha } \left(1+\frac{a^2}{x^2}+s^2\right)^{\inv{2}(p-\frac{3}{2})} e^{-|x|[\beta\left(1+\frac{a^2}{x^2}+s^2\right)^{\inv{2}}+\alpha\theta s]}\,ds \label{integral: integration-by-parts-GH}.
\end{align}

Hence, from (\ref{fxstar}) and (\ref{integral: integration-by-parts-GH}), as $ x \to - \infty$:
\begin{align}
\notag f_{X{*}}(x)=& \frac{e^{-\theta |x|}|x|^{p-\inv{2}}}{\sqrt{2\pi}\,\overline{K}_p(a, b) \beta^{p-\inv{2}}} 
\left(1+O\left(\inv{|x|}\right)\right) \\
& \quad \times \int^{\infty}_{\alpha } \left(1+\frac{a^2}{x^2}+s^2\right)^{\inv{2}(p-\frac{3}{2})} e^{-|x|[\beta\left(1+\frac{a^2}{x^2}+s^2\right)^{\inv{2}}+\alpha\theta s]}\,ds \label{fxstar1}
\end{align}

\section{Asymptotic bivariate equi-skew form}

We next need to investigate the asymptotic behaviour of the integral in (\ref{fxstar1}) as $x \to -\infty$.  To this end
define \begin{equation} \label{phistar}\phi(s) =    \beta\left(1+\frac{a^2}{x^2}+s^2\right)^{\inv{2}} + \alpha\theta s \end{equation}  where as before   $\beta = (\alpha^2\theta^2+\theta^2+b^2)^{\inv{2}},\alpha >0, \theta\in \R.$ 

We shall need in the sequel
\begin{equation} \label{index}
\phi(\alpha) + \theta >0, \, \, \phi(\alpha) >0,
\end{equation} To see this 
\begin{eqnarray*}
\phi(\alpha) - |\theta| &=& \beta\left(1+\frac{a^2}{x^2}+\alpha^2\right)^{\inv{2}} +\alpha^2\theta - |\theta|\\
&\geq& \beta(1+\alpha^2)^{\inv{2}} +\alpha^2\theta - |\theta|\\
&=& (1+\alpha^2)^{\inv{2}} ( (\alpha^2\theta+\theta^2+b^2)^{\inv{2}} +\alpha^2\theta- |\theta|\\
&>&(1+\alpha^2)^{\inv{2}} ( (\alpha^2\theta^2+\theta^2)^{\inv{2}} +\alpha^2\theta- |\theta|, \text{ since $b > 0$,}\\
&=& (1+\alpha^2)|\theta| + \alpha^2\theta - |\theta|\\
&=& \alpha^2(|\theta| + \theta) \geq 0.
\end{eqnarray*}
Hence $\phi(\alpha) > |\theta| > 0$, and (\ref{index}) follows. 
Next 
\begin{equation} \label{phistar1}
\phi^{'}(s) =\frac{\beta s}{\left(1+\frac{a^2}{x^2}+s^2\right)^{\inv{2}}}+\alpha\theta >0
\end{equation}
for all $s\geq \alpha$ with $\alpha$, $\beta>0$ and $\theta\in\R$ providing $|x| \geq \frac{|\theta a|}{b}$. To see this  we can show, similarly to the above, that
\begin{equation*}
\left(\beta s  - \alpha|\theta|\left(1+\frac{a^2}{x^2}+s^2\right)^{\inv{2}}\right) >0
\end{equation*} providing $s^2 > (\frac{\alpha \theta a}{xb})^2, $ so that,  if we take $s \geq \alpha$, providing $|x|>  \frac{|\theta a|}{b}.$
Given that we shall need $|x|\rightarrow \infty$, the inequality (\ref{phistar1})  will hold for all $s \geq \alpha$ for any fixed $\alpha>0, b >0, \theta \in \R.$

Thus in view of (\ref{phistar}), (\ref{phistar1}), $\phi(s)$, $s \geq \alpha,$ has  an inverse function $ \phi^{-1}(s), s \geq \phi(\alpha)$.
Next we consider,with reference to (\ref{fxstar1}),
\begin{equation} \label{seven}
  |x|^{p-\inv{2}}    \int^{\infty}_{\alpha } \left(1+\frac{a^2}{x^2}+s^2\right)^{\inv{2}(p-\frac{3}{2})} e^{-|x|\phi(s)}\,ds \end{equation} and change  variable of integration to $w = \phi(s)-\phi(\alpha)$, 
 so  the expression becomes:
\begin{eqnarray}&= &|x|^{p-\inv{2}} e^{-|x|\phi(\alpha)} \int^{w = \infty}_{w = 0} \frac{\left(1+\frac{a^2}{x^2}+ \left(\phi^{-1}\left(w+\phi(\alpha)\right)\right)^2\right)^{\inv{2}\left(p-\frac{3}{2}\right)}}{\phi'\left(\phi^{-1}(w+\phi(\alpha)\right)} e^{-|x|w}\,dw  \nonumber\\
&=& |x|^{p-\frac{3}{2}} e^{-|x|\phi(\alpha)}\{|x| \int^{w = \infty}_{w = 0} \tilde{\theta}(w)  e^{-|x|w}\,dw\} \label{thetatilde} \end{eqnarray}
where 
\begin{equation}
\tilde{\theta}(w) =  \frac{\left(1+\frac{a^2}{x^2}+\left(\phi^{-1}\left(w+\phi(\alpha)\right)\right)^2\right)^{\inv{2}\left(p-\frac{3}{2}\right)}}{\phi'\left(\phi^{-1}(w+\phi(\alpha)\right)}.  \label{tilde-theta1}
\end{equation}

We now consider each of the multiplicands in (\ref{thetatilde}) separately. First, using integration by parts, and putting for convenience $ v=|x|$ we have 
\begin{align*}
& v \int^{w=\infty}_{w=0} \tilde{\theta}(w)e^{-vw}\,dw \\
= &v \left[ \frac{\tilde{\theta}(w)e^{-vw}}{-v}\right]^{\infty}_{w=0} - v \int^{\infty}_{w=0} \frac{e^{-vw}}{-v} \tilde{\theta}'(w)\,dw \\
= & \tilde{\theta}(0) + \int^{\infty}_{w=0} \tilde{\theta}'(w)e^{-vw}\,dw.
\end{align*} 

We now assume for the moment  that $a=0$, so that  $\tilde{\theta}(w)$  does not involve   $v=|x|$; and  note in passing  that the case of the GH of special interest to us, the VG, would still be encompassed by an  initial assumption that $a=0$.

 From the the fact that $ \tilde{\theta}'(w) $ is bounded  near $w = 0^+$ and $ \tilde{\theta}'(w) $ for large positive $w$ is of asymptotic  growth: $\text{Const.}\times w^{\kappa} ,$ for some fixed $\kappa$ as $ w \to \infty$,  we obtain as $v \to \infty$:
\begin{equation*}
 \int^{\infty}_{w=0} \tilde{\theta}'(w)e^{-vw}\,dw= O\left(\frac{1}{v}\right)
\end{equation*} so that 
\begin{equation*} \label{thetalimit}  
 v \int^{w=\infty}_{w=0} \tilde{\theta}(w)e^{-vw}\,dw =   \tilde{\theta}(0)\left(1+O\left(\inv{v}\right)\right), \,  v \to \infty,
\end{equation*} where  from (\ref{tilde-theta1}) and (\ref{phistar1}}) with $a=0$
\begin{align}
\tilde{\theta}(0) & = \frac{\left(1+\alpha^2\right)^{\inv{2}(p-\frac{3}{2})}}{\phi'(\alpha)} = \frac{(1+\alpha^2)^{\inv{2}(p-\frac{3}{2})}}{\frac{\beta\alpha}{(1+\alpha^2)^{\inv{2}}}+\alpha \theta}  = \frac{(1+\alpha^2)^{\inv{2}(p-\inv{2})}}{\alpha(\beta+\theta(1+\alpha^2)^{\inv{2}})}. \label{theta0}
\end{align}

 Thus from  (\ref{seven}) 
\begin{equation} \label{seven1}
  |x|^{p-\inv{2}}    \int^{\infty}_{\alpha } \left(1+\frac{a^2}{x^2}+s^2\right)^{\inv{2}(p-\frac{3}{2})} e^{-|x|\phi(s)}\,ds =\tilde{\theta}(0) |x|^{p-\frac{3}{2}}  e^{-|x|\phi(\alpha)} \left(1+O\left(\inv{|x|}\right)\right) \end{equation} when $a=0$ as $ x\to - \infty$, where
$\tilde{\theta}(0)$ is given by (\ref{theta0}), and 
\begin{equation} \label{phialpha}
\phi(\alpha) = \beta\left(1+\alpha^2\right)^{\inv{2}} + \alpha^2\theta 
\end{equation}

 It can be shown that (\ref{seven1}) holds with these same values of  $\tilde{\theta}(0)$ 
and $\phi(\alpha)$ {\it for general}  $a$.
Hence, for GH, 
\begin{equation*}
 f_{X^{*}}(x) =\frac{\tilde{\theta}(0) e^{-\theta |x|}|x|^{p-\frac{3}{2}}  e^{-|x|\phi(\alpha)}}{\sqrt{2\pi}\,\overline{K}_p(a, b) \beta^{p-\inv{2}}}\left(1+O\left(\inv{|x|}\right)\right) \label{fxstar3}
\end{equation*}
as $ x\to - \infty$.\\
Thus for large negative $y$, noting  from (\ref{index})   that $ \phi(\alpha)+ \theta >0$, so the integral is well-defined:
\begin{align*}
P (X^{*} \leq y)& = \frac{\tilde{\theta}(0)}{\sqrt{2\pi}\,\overline{K}_p(a, b) \beta^{p-\inv{2}}}
\int_{-\infty}^ y |x|^{p-\frac{3}{2}}  e^{-|x|(\phi(\alpha)+ \theta)}\left(1+O\left(\inv{|x|}\right)\right) dx \nonumber \\
&= \frac{\tilde{\theta}(0)}{\sqrt{2\pi}\,\overline{K}_p(a, b) \beta^{p-\inv{2}}}
\int_{|y|}^ {\infty} v^{p-\frac{3}{2}}  e^{-v(\phi(\alpha)+ \theta)}\left(1+O\left(\inv{v}\right)\right) dv \label{xstar4}
\end{align*}
where $v = -x = |x|$. 

Now, from L'H\^{o}pital's rule,  for $ \alpha  \in \R,  \beta >0 $
\begin{eqnarray} 
\int_s^\infty y^{\alpha - 1} e^{-\beta y} dy \, &=& \,  \frac{1}{\beta} s^{\alpha - 1} e^{-\beta s}\left(1 + o(1)\right), \,  s \to \infty,  \label{thesis.p26}  \\
&=&\frac{1}{\beta} s^{\alpha - 1} e^{-\beta s}+ \frac{\alpha - 1}{\beta} \int_s^\infty y^{\alpha - 2} e^{-\beta y} dy \nonumber
\end{eqnarray}  by integration by parts. Thus
\begin{eqnarray*} 
\left|\int_s^\infty y^{\alpha - 1} e^{-\beta y} dy -  \frac{1}{\beta} s^{\alpha - 1} e^{-\beta s}\right| 
&\leq& \frac{|\alpha - 1|}{\beta s} \int_s^\infty y^{\alpha - 1} e^{-\beta y} dy\\
&=& \frac{|\alpha - 1|}{\beta s}  \frac{1}{\beta} s^{\alpha - 1} e^{-\beta s}\left(1 + o(1)\right),\quad  \text{from (\ref{thesis.p26}), }\\
&=& \frac{1}{\beta} s^{\alpha - 1} e^{-\beta s}\left(O\left(\frac{1}{s}\right)\right)
\end{eqnarray*}   so that 
\begin{equation*}\label{thesisextn}
\int_s^\infty y^{\alpha - 1} e^{-\beta y} dy =  \frac{1}{\beta} s^{\alpha - 1} e^{-\beta s}\left(1+O\left(\frac{1}{s}\right)\right), \, s \to \infty,
\end{equation*} whence ,  as $ y \to -\infty$:
\begin{equation} \label{xstar5}
P (X^{*} \leq y) =  \frac{\tilde{\theta}(0)}{\sqrt{2\pi}\,\overline{K}_p(a, b) \beta^{p-\inv{2}}( \phi(\alpha)+ \theta) } |y|^{p-\frac{3}{2}}  e^{-|y|(\phi(\alpha)+ \theta)}\left(1+O\left(\inv{|y|}\right)\right) 
\end{equation} where
$\tilde{\theta}(0)$ is given by (\ref{theta0}), and $\phi(\alpha) + \theta =(1+\alpha^2)^{\inv{2}}(\beta+\theta(1+\alpha^2)^{\inv{2}}),$ and recall:
\begin{equation} \label{equality}
P(X_1 \leq y, \, X_2 \leq y)=P (X^{*} \leq y).
\end{equation}
\section{Quantile function and asymptotic copula}
The marginal density of each of $X_1, X_2$ is given by
\begin{equation*}
f_{X_1}(x) = \frac{e^{\theta x}}{\sqrt{2\pi}\, {\overline{K}}_p(a,b)}\overline{K}_{p-1/2}((x^2 +a^2)^{1/2},(\theta^2 +b^2)^{1/2}), x \in  \R,
\end{equation*} as expressed in \cite{Fung2011a}, equation (20), following \cite{Blaesild1981}. Hence after some algebra using (\ref{bessel})
\begin{equation} \label{ghcdf}
F_1(x) = P(X_1\leq x) = A |x|^{p-1}e^{-\left((\theta^2+b^2)^{\inv{2}}+\theta\right)|x|}\left(1+O\left(\inv{|x|}\right)\right)
\end{equation}
where  $A^{-1}= 2\overline{K}_p(a,b)\left(\theta^2+b^2\right)^{\frac{p}{2}} \left((\theta^2+b^2)^{\inv{2}}+\theta\right)$, as $x\to -\infty$. In the VG special case where $a=0, b =\sqrt{\frac{2}{\nu}}, p=1/2, \mu=0, \sigma^2=1,$ this is equation (18) of \cite{Fung2011a}.

The expression (\ref{ghcdf}) itself is a special case of distribution functions with a generalized gamma-type tail considered by \cite{FungSeneta2018}, equation (6), where in the notation of that paper on the left-hand  side of the following:
$a=A, b=p-1, c = (\theta^2 + b^2)^{\inv{2}}+ \theta, d=1, e=1.$ So from (9) of  that paper, in our notation, we have 
\begin{equation*}
F_1^{-1}(u)  = \frac{\log u}{(\theta^2+b^2)^{\inv{2}}+\theta} - \frac{\left(p-1\right)\log|\log u|}{(\theta^2+b^2)^{\inv{2}}+\theta} - \frac{\left(p-1\right)\log\left(\frac{A^{\frac{1}{p-1}}}{(\theta^2+b^2)^{\inv{2}}+\theta}\right)}{(\theta^2+b^2)^{\inv{2}}+\theta}+ O\left(\frac{\log|\log u|}{|\log u|}\right)
\end{equation*}
as $u\to 0^+$.

We finally address the rate of convergence.  To simplify notation put 
\begin{align*}
\gamma &= \left(\beta + (1+\alpha^2)^{\inv{2}}\theta\right) \quad \left(= \left(\phi(\alpha)+\theta\right)/\left(1+\alpha^2\right)^{\inv{2}}\right)\\
\delta &=(\theta^2+b^2)^{\inv{2}}+\theta\\
\tau &= \frac{(1+\alpha^2)^{\inv{2}}\gamma}{\delta}\\
C_1 &= \frac{(1+\alpha^2)^{\inv{2}\left(p-\frac{3}{2}\right)}}{\sqrt{2\pi}\overline{K}_{p}(a,b)\beta^{p-\inv{2}}\alpha(1+\alpha^2)^{\inv{2}}\gamma^2\delta^{p-\frac{3}{2}}}\\
C_2 &= A^{- \tau}\delta^{(p-1)\tau}
\end{align*}

Then from  (\ref{xstar5}), (\ref{equality}) we have {\adb
\begin{eqnarray*}
&& P(X_1\leq F_1^{-1}(u),\, X_2\leq F_2^{-1}(u)) \\
&\sim& C_1\left(\left| \log u + O\left(\log|\log u|\right)\right| \right)^{p-\frac{3}{2}}
e^{- \left|\frac{\log u}{\delta} - \frac{\left(p-1\right)\log|\log u|}{\delta} - \frac{\left(p-1\right)\log C_2}{\delta} + O\left(\frac{\log|\log u|}{|\log u|}\right)\right| (1+\alpha^2)^{\inv{2}}\gamma}\\
&\sim& 
C_1 \times u^{\tau} \times 
\left(|\log u|\right)^{\left(p-\frac{3}{2}\right)-\left(p-1\right)\tau}\times C_2.
 \end{eqnarray*}

As a result, 
\begin{equation*}
 \frac{P(X_1\leq F_1^{-1}(u),\, X_2\leq F_2^{-1}(u))}{u}  =  u^{\tau - 1} L(u)
\end{equation*}
where 
\begin{equation*}
L(u) \sim  C_1C_2 \left(|\log u|\right)^{\left(p-1\right)\left(1-\tau \right)- \inv{2}}
\end{equation*}
is a slowly varying function.

Obviously,  the rate for the VG is obtained explicitly by letting $a=0$, $b = \sqrt{\frac{2}{\nu}}$, $p=\inv{2}$.

\bibliography{library}{}

\end{document}